\documentclass{article}

\usepackage{amsfonts,amsmath,amssymb,amsthm}
\usepackage[all]{xypic}

\newcommand{\Cc}{\mathbb{C}}  

\newtheorem*{theorem*}{Theorem}  
\newtheorem{theorem}{Theorem}    
\newtheorem{lemma}{Lemma}        

\begin{document}

\title{\textbf{Milnor fibration and fibred links at infinity}}
\author{Arnaud Bodin}
\date{29 January 1999}
\maketitle

\section*{Introduction}

Let  $f : \Cc^2 \longrightarrow \Cc$ be a  polynomial function.
By definition  $c\in \Cc$
is a \emph{regular value at infinity} if there exists a disc $\mathcal{D}$
centred at $c$ and a compact set $\mathcal{C}$ of $\Cc^2$ 
such that the map
 $f : f^{-1}(\mathcal{D})\setminus \mathcal{C} \longrightarrow \mathcal{D}$ 
is a locally trivial fibration. There are only a finite number of 
\emph{critical} (or \emph{irregular}) \emph{values at infinity}.
For $c \in \Cc$ and  a sufficiently large real number $R$, the
 \emph{link at infinity}  $K_c = f^{-1}(c) \cap S^3_R$ is well-defined.

In this paper we sketch the proof of the  following
theorem which gives a characterization of  fibred multilinks 
at infinity.
\begin{theorem*}
A multilink $K_0 = f^{-1}(0)\cap S^3_R$ is fibred 
if and only if all the values $c \not= 0$ are regular at infinity.
\end{theorem*}

We first obtain theorem \ref{th:fib}, a version of this theorem 
 was proved by A.~N\'emethi and A.~Zaharia in \cite{NZ} 
(with ``semitame'' as a hypothesis). Here we give 
a new proof using   resolution  of 
singularities at infinity. This method enables us  to describe the
fibre and the  monodromy of the Milnor fibration in terms of
combinatorial invariants of a resolution of $f$.

\begin{theorem}
\label{th:fib}
If there is no critical value at infinity outside $c=0$ then in the homotopy
class of
$$
 \frac{f}{|f|} : S^3_R \setminus f^{-1}(0) \longrightarrow S^1 
$$
there exists a fibration.
\end{theorem}

The value $0$ may be regular or not. One may specify what kind of fibration it is; 
if $f$ is a reduced polynomial, then this is an open book decomposition, 
otherwise it is a multilink fibration of $K_0 = f^{-1}(0) \cap S^3_R$ 
(see paragraph \ref{sec:def}). 
The weights of $K_0$ are given by the multiplicities of the factorial
decomposition of $f$.

If $0$ is a regular value at infinity and $c \not=0$ is a critical value 
at infinity, W.~Neumann and L.~Rudolph proved in \cite{NR} that the link $f^{-1}(0)\cap S^3_R$ 
is not fibred. In the following theorem~\ref{th:nonfib} we do not have any
hypothesis on the value $0$, in particular $0$ can be a critical 
value at infinity.

\begin{theorem}
\label{th:nonfib}
Suppose that $c\not= 0$ is a critical value at infinity for $f$, then
 the multilink  $K_0 = f^{-1}(0) \cap S^3_R$
is not a fibred multilink.
\end{theorem}

We begin with definitions, the second part is devoted to the proof of 
theorem~\ref{th:fib}. We conclude with the proof of theorem~\ref{th:nonfib}.

\section{Definitions}
\label{sec:def}

As in \cite{EN}, a \emph{multilink}  $L(\mathbf{m})$ ($\mathbf{m} = (m_1,\ldots,m_k)$) is
 a link having each component $L_i$ weighted by the integer $m_i$.

The multilink  $L(\mathbf{m})$ is a \emph{fibred multilink} if
there exists a differentiable fibration 
$\theta : S^3_R\setminus L \longrightarrow S^1$
such that $m_i$ is the degree of the restriction of
$\theta$ on a meridian of $L_i$. A fibre
$\theta^{-1}(x)$ is a \emph{Seifert surface} for the multilink.
The link $K_0 = f^{-1}(0) \cap S^3_R$ is  a multilink, the weights being naturally
given by the multiplicities of the factorial decomposition of $f$.

A \emph{fibred link} is a fibred  multilink
having all its components weighted by $+1$. Then  $\theta$ is called an
\emph{open book decomposition}.

\medskip

Next we give definitions and results about resolutions, see \cite{LW}.
Let $n$ be the degree of $f$ and $F$ be the corresponding homogeneous polynomial 
with the same degree.
The map $\tilde{f} : \Cc P^2 \longrightarrow \Cc P^1$, $\tilde{f}(x:y:z) = (F(x,y,z):z^n)$ 
is not everywhere defined, nevertheless there exists a minimal composition
of blowing-ups  $\pi_w : \Sigma_w \longrightarrow \Cc P^2$
such that  $\tilde{f} \circ \pi_w$ extends to a well-defined morphism
$\phi_w$  from  $\Sigma_w$ to  $\Cc P^1$.
This is the \emph{weak resolution}.
$$
\xymatrix{
{}\Cc^2   \ar[r] \ar[d]_f       &\Cc P^2 \ar[d]_{\tilde{f}} 
    &\Sigma_w \ar[l]_-{\pi_w} \ar[dl]^{\phi_w}  \\
\Cc \ar[r]      &\Cc P^1
}
$$

\bigskip

For an irreducible component $D$ of $\pi_w^{-1}(L_\infty)$
 ($L_\infty$ is the line of $\Cc P^2$ having the 
equation $(z=0)$), we distinguish three cases: 
\begin{enumerate}
  \item $\phi_w(D)=\infty$, we denote  $D_\infty = \phi_w^{-1}(\infty)$.
  \item $\phi_w(D) = \Cc P^1$, $D$ is a \emph{dicritical component},
 the restriction of
 $\phi_w$ to $D$ is a ramified covering,
the \emph{degree} of $D$ is the degree of this restriction. 
The divisor which contains all
these components is the \emph{dicritical divisor} 
$D_{\!\text{\textit{dic}}}$.
  \item $\phi_w(D) = c \in \Cc$, there is a finite number of such components,
collected in  $D_{\!\text{\textit{crit}}} = D_{c_1}\cup\ldots\cup D_{c_g}$.
\end{enumerate}

The irregular values at infinity for $f$ are the values  $c_1,\ldots,c_g$
and the critical values of the map $\phi_w$ restricted to  
$D_{\!\text{\textit{dic}}}$;
moreover each divisor $D_{c_i}$ is a disjoint union of bamboos.

\medskip

We now increase the number of blowing-ups of $\pi_w$ in a minimal way, in order 
to obtain $\pi_p : \Sigma_p \longrightarrow \Cc P^2$ and
$\phi_p = \tilde{f} \circ {\pi_p} : \Sigma_p \longrightarrow \Cc P^1$ such
that the fibre $\phi_p^{-1}(0)$ cuts the divisor $D_{\!\text{\textit{dic}}}$ 
transversally and is a normal crossing divisor. This is the 
\emph{partial resolution} for the value $c=0$.

We continue with blowing-ups in order to obtain $\pi_t, \Sigma_t, \phi_t$ such
that each fibre of $\phi_t$ cuts the divisor $D_{\!\text{\textit{dic}}}$ 
transversally and all the fibres of $\phi_t$
are normal crossing divisors.
This is the \emph{total resolution}.

For the total resolution the values  $c_1,\ldots,c_{g'}$ coming from the 
components $D$ of the new $D_{\!\text{\textit{crit}}}$ with $\phi_t(D) = c_i$
are the critical values at infinity.

\section{Milnor fibration at infinity}

Until the end of this section, we suppose that the only irregular
value at infinity for $f$ can be the value $0$. Let $\phi = \phi_t$ coming 
from the total resolution.
In $\Sigma_t$ the sphere $\pi_t^{-1}(S^3_R)$  is diffeomorphic to the boundary $S$ of
a neighbourhood of $\pi_t^{-1}(L_\infty)$ (see \cite{D}).

Instead of studying $f/|f|$ restricted to $S^3_R \setminus f^{-1}(0)$
we study $\phi/|\phi|$ restricted to $S \setminus \phi^{-1}(0)$.
Let $\theta$ be the restriction of $\phi/|\phi|$ to $S \setminus
 \phi^{-1}(0)$. By changing the sphere $\pi_t^{-1}(S^3_R)$ into $S$  we only know
that $\theta$ is in the homotopy class of $f/|f|$.

As in \cite{LMW} there is a correspondence between the irreducible components
of $\pi_t^{-1}(L_\infty)$ and a Waldhausen decomposition of $S\setminus \phi^{-1}(0)$
into Seifert three-manifolds. We will prove that the restriction of $\theta$ to the 
Seifert manifold $\sigma(D)$ associated to any irreducible component $D$ of 
$\pi_t^{-1}(L_\infty)$ is a fibration. If $D \subset D_\infty \cup D_0$, the
equations are similar to the local case; we thus have to look at what happens
with the components of the dicritical divisor.

\begin{lemma}
\label{lem:ann}
The smooth points in $\pi_t^{-1}(L_\infty)$ of each dicritical component with non-empty 
intersection with $D_{\!\text{\textit{crit}}}=D_0$ is an annulus.
\end{lemma}
In other words the intersection of $D_0$ with each dicritical
component is empty or reduced to a single point.

\begin{proof}
 This is a consequence of the fact that above
$\Cc P^1 \setminus \{0,\infty\}$, $\phi$ is a regular covering.
\end{proof}

With similar arguments, one can prove:
\begin{lemma} 
\label{lem:disc}
Each dicritical component $D$ with $D \cap D_{\!\text{\textit{crit}}} = \varnothing$ is
of degree $1$.
\end{lemma}

\subsection{Fibration on $\sigma(D)$ for $D \subset D_{\!\text{\textit{dic}}}$ }

Let $D$ be a dicritical component and let $U$ be the simple points of $D$
in $\pi_t^{-1}(L_\infty) \cup \phi^{-1}(0)$. By
lemmas \ref{lem:ann} and \ref{lem:disc} 
we know that $U$  is an annulus and 
 $\phi_{|{U}} : {U}  \longrightarrow \Cc P^1\setminus \{0,\infty\}$
is a regular covering of order $d$.

Let $u \in \Cc^*$ be a parametrisation of $U$. For each point of $U$ we choose
local coordinates $(u,v)$ such that
$\phi$ can be written $\phi(u,v) = u$. We choose $S$ so that
$S$ is locally given by $(|v|=\varepsilon)$
where $\varepsilon$ is a small positive real number.

With these facts one can calculate that 
the restriction of $\theta$ to the Seifert component $\sigma(D)$ associated to $D$ 
is a fibration whose fibres consist of $d$ annuli.

\subsection{Fibration in a neighbourhood of a non-simple point}

In a neighbourhood $V$ of a non-simple point, i.e.~a point
belonging to  a dicritical component $D$ and another component 
$D' \in \pi_t^{-1}(L_\infty) \cup \phi^{-1}(0)$, $\phi$ is
defined in appropriate local coordinates by  $(u,v) \mapsto u^d$.

Let $T$ be the tubular neighbourhood of $D\cap V$ given by $(|v|\leqslant \varepsilon)$.
$\theta_{|T}$  defines a fibration whose fibres consist of $d$ annuli:
$$
 \theta^{-1}(e^{i\alpha}) \cap T = 
 \Big\lbrace (u,v) \in T; |v| = \varepsilon, u \not= 0 \text{\ \ and\ \ }
 {u^d}/{|u|^d}=e^{i\alpha} \Big\rbrace.
$$
 
For $T'$ a tubular neighbourhood of $D'\cap V$ given by $(|u| \leqslant \varepsilon)$, 
the fibre $ \theta^{-1}(e^{i\alpha}) \cap T'$ is also a union of $d$ annuli.

These different pieces fit nicely on the torus  $\partial T \cap \partial T'$. 
So  with a plumbing of $T$ and $T'$, $\theta$ is a fibration on $V$.

\subsection{Fibration in a neighbourhood of the strict transform}

Let $F$ be an irreducible component of $\phi^{-1}(0)\setminus D_0$ (which
corresponds to the affine set $f^{-1}(0)$). $F$ can intersect $D_0$ or $D_{\!\text{\textit{dic}}}$.
If $F\cap D_0 \not= \varnothing$ then locally in a neighbourhood $V$, 
$\phi(u,v)=u^pv^q$ with $(v=0)$ is an equation for $D_0$. The associated 
component of the link is 
$\phi^{-1}(0)\cap S \cap V$. Then  $\theta_{|V}$ is a fibration whose fibres 
consist of $\mathrm{gcd}(p,q)$ annuli:  
$$
 \theta^{-1}(e^{i\alpha}) \cap V = 
 \Big\lbrace (u,v) \in V; |v| = \varepsilon, u \not= 0 \text{\ \ and\ \ } 
 {u^pv^q}/{|u^pv^q|}=e^{i\alpha} \Big\rbrace.
$$ 
Moreover this fibration is a multilink fibration, because on a torus
 $D^2_\delta \times S^1_\varepsilon \setminus\{0\}$, 
the trace of the fibre at $v=\text{\textit{cst}}$ is $p$ radii of the annulus
$D^2_\delta \setminus \{ 0 \} \times v$.
If $f$ is a reduced polynomial function then $p=1$ and $\theta$ is an open book decomposition.

Similarly, $\theta$ is still locally a fibration  if 
$F\cap D_{\!\text{\textit{dic}}} \not= \varnothing$.

\medskip

We now conclude by collecting  and gluing previous results. 
${\phi}/{|\phi|}$ is a fibration in a neighbourhood of $S\cap \phi^{-1}(0)$
and on all $V \cap S$ which cover  $S \setminus \phi^{-1}(0)$, so 
${\phi}/{|\phi|} : S \setminus \phi^{-1}(0) \longrightarrow S^1$
is a fibration.
Furthermore with the discussion above ${\phi}/{|\phi|}$
is an open book decomposition or a multilink fibration
depending on  $f$ being reduced or not.

\section{Non-fibred multilinks}

Under the hypotheses of theorem \ref{th:nonfib} and without loss of generality 
we suppose that $\{\lambda c \text{ with } \lambda < 0 \}$ 
does not contain critical values of $f$ at infinity.
The surface 
 $\mathcal{F} = \left( {f}/{|f|}\right)^{-1} (-{c}/{|c|}) \cap S^3_R$ 
is a Seifert surface for the multilink $K_0=f^{-1}(0)\cap S^3_R$.
Moreover, for complex numbers $\omega$ with $0 \leqslant |\omega-c| \ll |c|$
the links  $f^{-1}(\omega)\cap S^3_R$ do not cut $\mathcal{F}$.

We choose $\omega$ as a regular value at infinity. 
For the partial resolution $\phi=\phi_p$ at infinity for $f$ and the value $0$, 
there exists one dicritical
component with a valency at least $3$ in $\pi_p^{-1}(L_\infty) \cup \phi^{-1}(0)$:
let $D$ be a dicritical component where $c$ is a critical value at infinity.
If the intersection $\phi^{-1}(0) \cap D$ has more than two points or if there is
a bamboo of $D_c$ that cuts $D$ then we can easily conclude. But no other case is 
possible because
$\phi_{|D}$, with the critical values $0$ and $c$, has more than two zeroes.
So the manifold $\sigma(D)$ induces a Seifert manifold of the minimal
decomposition of $S^3_R \setminus K_0$;
by crossing this component, $f^{-1}(\omega)$ 
defines a \emph{virtual component} of  $M = S^3_R \setminus f^{-1}(0)$ (see \cite{LMW}):
that is to say a regular fibre of the minimal 
Waldhausen decomposition of the manifold $M$.

According to \cite[th.~11.2]{EN}, since  $\mathcal{F}$ 
and a virtual component of $M$ have empty intersection,
 $K_0$ is not a fibred multilink, if we exclude the case where $M$
is $S^1\times S^1 \times [0,1]$. This case is studied in the
following lemma.

\begin{lemma}
If the underlying link associated to $K_0=f^{-1}(0)\cap S^3_R$ is the
Hopf link then $c \not= 0$ is a regular value at infinity for $f$.
\end{lemma}

\begin{proof}
We suppose first that $f$ is a
reduced polynomial function. Then $K_0$ is the Hopf link, and since $K_0$ is an 
iterated torus link around Neumann's multilink $L$ \cite[\S2]{N}, 
this multilink can only be the trivial knot or the Hopf
link.

\emph{Case of $L$ being the trivial knot:} 
There is only one dicritical component. 
If $f$ is not a primitive polynomial (i.e.~with connected generic fibre) 
then with the use of the Stein factorisation, let  $h \in \Cc[t]$ and let
 $g \in \Cc[x,y]$ be a primitive polynomial  with  $f = h \circ g$.
By the  Abhyankar-Moh theorem (see \cite{A}), there exists an algebraic automorphism
$\Theta$ of  $\Cc^2$ with  $g\circ \Theta (x,y)=x$ and then
 $f \circ \Theta(x,y)= h(x)$.

Let $x_1,\ldots,x_n$ be the zeroes of $h$; 
$x_1\times\Cc,\ldots,x_n\times\Cc$ are the solutions of 
$f\circ \Theta (x,y)=0$. Therefore the link $K_0$ is a union
of trivial knots with zero linking numbers, so $K_0$ is 
not the Hopf link.

\emph{Case of $L$ being the Hopf link:}
 $K_0$ and the multilink $L$ are isotopic. 
On the one hand in the weak resolution for $f$, the restriction of $\phi = \phi_w$
to $D_{\!\text{\textit{dic}}}$
cannot  have the critical value $0$ without a bamboo. If so, one component of $K_0$ would
be a true iterated torus knot around a component of $L$, in contradiction
with the hypothesis. On the other hand, 
each component of the multilink $L$ can be represented by a disc which crosses transversally
the last component of each bamboo (start counting at the dicritical component).
If there exists a bamboo for the value $0$, the component $C$ of 
$\phi^{-1}(0)\setminus D_0$ with $C\cap D_0 \not= \varnothing$
must be irreducible, reduced and cross $D_0$ transversally at the last component; this configuration
is excluded by lemma 8.24 of \cite{MW}. So $0$ is a regular value
at infinity and since $K_0$ is isotopic to $L$, all the dicritical
components have degree one and there is no value having a bamboo, so
$c$ is a regular value at infinity for $f$.

If $f$ is not reduced, let $g$ be the reduced polynomial function associated
to $f$. Then the link $g^{-1}(0)\cap S^3_R$ is the Hopf link and from
the discussion above we know that $0$ is a regular value at infinity for $g$.
From the classification of regular algebraic annuli \cite[\S8]{N},
there exists an algebraic automorphism $\Theta$ of $\Cc^2$ with 
$\Theta(0,0) = (0,0)$ such that
$g \circ \Theta (x,y) = xy+\lambda$, $\lambda \in \Cc$. So $f\circ \Theta (x,y) = (xy+\lambda)^l$ if
$\lambda \not= 0$ and $f\circ \Theta (x,y) = x^py^q$ if $\lambda = 0$. In both cases,
$c$ is a regular value at infinity for $f$.
\end{proof}

In conclusion, whether $0$ is a regular value at infinity 
or not, the multilink $K_0 = f^{-1}(0) \cap S^3_R$ is not fibred
when $c \not= 0$ is a critical value at infinity.

\bigskip

\section*{Acknowledgments}

{\small
I thank Professor~Fran\c{c}oise~Michel for long discussions 
and for her many ideas.
}

{\small
Universit\'e Paul Sabatier Toulouse III, laboratoire \'Emile Picard, 
118 route de Narbonne, 31062 Toulouse cedex 4, France;
bodin@picard.ups-tlse.fr
}

\end{document}